\documentclass[a4paper, 10pt,conference]{IEEEtran}
\usepackage[lining]{ebgaramond}
\IEEEoverridecommandlockouts
\usepackage{subfigure}
\usepackage[margin=1in]{geometry}
\usepackage[utf8]{inputenc}
\usepackage{cite}
\usepackage{mwe}
\usepackage{amsmath}
\usepackage{mathrsfs}
\usepackage{url}
\usepackage{microtype}
\usepackage{dirtytalk}
\usepackage{eucal}
\usepackage{amsmath,amssymb,amsfonts}
\usepackage{algorithmic}
\usepackage{graphicx}
\usepackage{textcomp}
\usepackage{xcolor}
\usepackage{algorithm}
\usepackage{color}
\usepackage{multirow}
\usepackage{listings}
\usepackage{graphicx}
\usepackage{amsmath}
\usepackage{hyperref}
\usepackage{authblk}
\DeclareFontFamily{U}{calligra}{}
\DeclareFontShape{U}{calligra}{m}{n}{<->callig15}{}

\title{A closed loop gradient descent algorithm applied to Rosenbrock's function}
\author[$*$]{Subhransu S. Bhattacharjee\thanks{$*$Mr. Subhransu S. Bhattacharjee is a student at the Australian National University, Canberra, Australia. Please direct all queries to Subhransu Bhattacharjee at u7143478@anu.edu.au.}}
\author[$\dagger$]{Ian R. Petersen\thanks{$\dagger$Dr. Ian R. Petersen, FAA is a professor at the College of Engineering and Computer Science, Australian National University, Canberra, Australia.}}
\affil[$*\,\dagger$]{The Australian National University, Canberra, ACT 0200, Australia}
\begin{document}
\maketitle
\noindent
\begin{abstract}
We introduce a novel adaptive damping technique for an inertial gradient system which finds application as a gradient descent algorithm for unconstrained optimisation. In an example using the non-convex Rosenbrock's function, we show an improvement on existing momentum-based gradient optimisation methods. Also using Lyapunov stability analysis, we demonstrate the performance of the continuous-time version of the algorithm. Using numerical simulations, we consider the performance of its discrete-time counterpart obtained by using the symplectic Euler method of discretisation.
\end{abstract}
\noindent
\section{Introduction} \label{sec:introduction}
 \noindent Recent advances in the field of deep learning has rekindled interest in gradient-based optimisation algorithms. Often we find the training loss for machine learning models is dependent on the nature of such algorithms. This motivates us to study these algorithms from a control systems perspective. While studying such algorithms, we encounter unconstrained optimisation problems which are of the form:
\begin{equation}\label{1}
\mathop{min}_{x \in \mathbb{R}^d} f(x) \;\;\;\; \text{where}\, f: \mathbb{R}^{d}\longrightarrow{\mathbb{R}}\;.
\end{equation}
Mostly, such algorithms are dependent on learning the gradient of the cost function. A central aspect of analysing this problem is the Lipschitz continuity of the cost function \cite{hiji}, where there exists a constant $L$, termed as the Lipschitz constant, such that:
\begin{equation}\label{2}
  \| \nabla{f}(y)-\nabla{f}(x) \|\leq L \| y-x \|  
\end{equation}
$\forall \; x,y\,\in \mathbb{R^d} $. A sufficient condition for learning the gradient of such a cost function is to take a small enough step-size $s$ such that:
\begin{equation}\label{3}
    0<s\leq \frac{1}{L}\,.
\end{equation}
Thus, we consider methods which require first order gradient knowledge to obtain the minima. However, the performance of such methods is heavily dependent upon the cost function's spectral condition number, its overall geometry, the presence of saddle points, and local minimas \cite{bert}. \\[0.1cm]
One such method is the Nesterov's accelerated gradient descent algorithm which is given by:
\begin{equation}\label{4}
\begin{split}
    &     \theta_{k+1} = \theta_k + v_{k+1},\\
    &      v_{k+1} = \theta_{k+1} - \theta_k,\\
    & \theta_k = y_k - \varepsilon_k \nabla f\left(y_k\right),\\
    & y_{k+1} = \theta_k + \mu_k\left(\theta_k - \theta_{k-1} \right),
\end{split}
\end{equation}
with a reported convergence rate of $\mathcal{O}(\frac{1}{k^2})$ \cite{10029946121} for a convex cost function $f$. In continuous-time this algorithm takes the form [11]:
\begin{equation}\label{5}
    \ddot{X}+\frac{\alpha}{t}\,\dot{X}+\nabla f(x) = 0\;\;\; \forall\; t>0.
\end{equation}
\begin{figure}
\centering
\includegraphics[width=0.7\columnwidth]{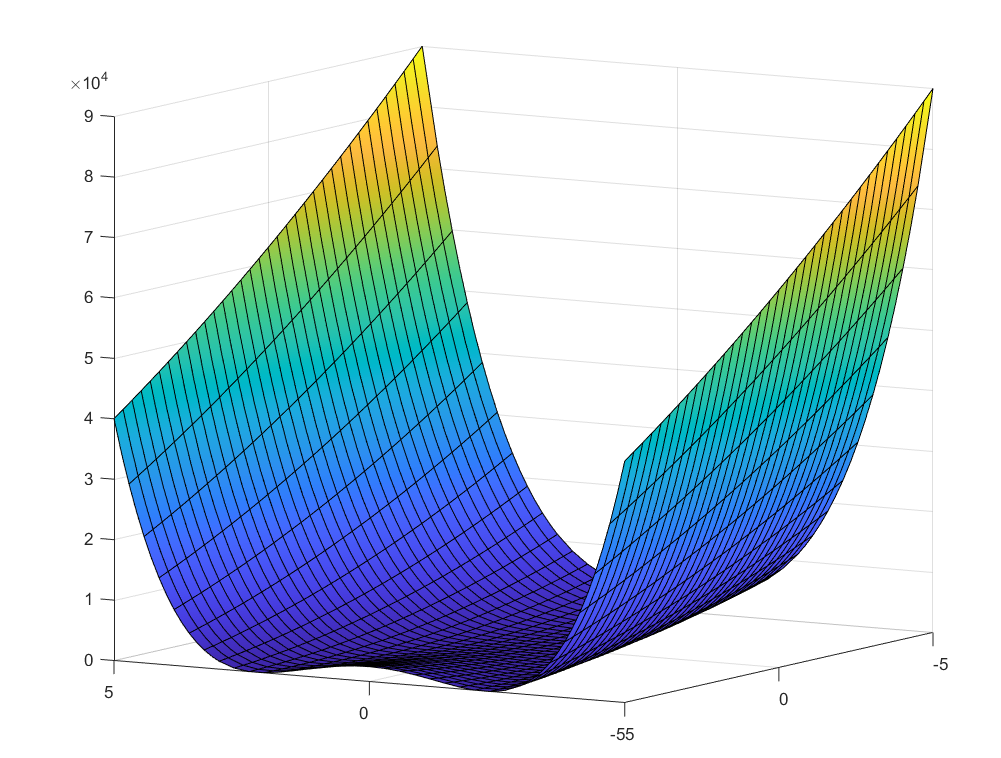}
 \caption{\small Rosenbrock's function with global minima at $(1,1)$.}
 \label{fig:RF}
\end{figure}
\noindent 
Su et al. show in their paper \cite{su2015differential} how the possible values of the damping coefficient, $\alpha$ affect the rate of convergence of this algorithm. For convex functions, and $\alpha \geq3$, the rate of convergence of (\ref{5}) is of the order $f-f^* = \mathcal{O}(\frac{1}{t^2})$. Though there is a lack of mathematical literature to explain the acceleration of this algorithm, recent advances show a variety in the rate of convergence by re-scaling the gradient flow which corresponds to a Bregman Lagrangian \cite{wibi}. While this does provide some insight, all derivations in the class of inertial gradient systems like (\ref{5}) are based on realising Lyapunov functionals which themselves rely on the geometry of the cost function to approximate an upper bound on the convergence rate. These algorithms in continuous-time can be regarded as open loop systems which show variance in their dynamic evolution in discrete-time, depending upon the order of discretization used. The Nesterov scheme for example takes multiple forms depending upon the discretisation method used \cite{muehlebach2019dynamical}. \\
\noindent
Notwithstanding the difficulty in explaining the phenomenon of acceleration in the Nesterov and Heavy-ball \cite{aujol:hal-02545245} schemes, a number of studies \cite{fist} have shown that momentum definitely plays a central role in the acceleration of such optimisation algorithms. It should be noted that the continuous-time analysis of such systems is achieved by using high-resolution ordinary differential equation approximations with small step sizes. The literature on first order momentum-based methods shows that, in the deterministic setting, such methods reveal a stabilising effect in their transient phases. Recent results indicate that these momentum methods admit an attractive invariant manifold on which the dynamics reduce to a gradient flow \cite{kovachki2021continuous}.
\\[0.1cm]
In this paper, we will discuss a particular method within a relatively new class of gradient descent schemes which are closed loop in nature and show qualitatively better numerical performance in continuous-time compared to an open loop approach like Nesterov's scheme. Our study begins with understanding the problem generalised as:
\begin{equation}\label{7}
    \ddot{X}\,+\,\gamma\,\dot{X}+\,\nabla f\,(x)\,\,=\,0,
\end{equation}
where $\gamma$ is the damping coefficient. We will consider the construction of $\gamma$ as a feedback control problem to optimise the dynamics for rapid convergence.\\[0.1cm]
To test the performance of our method, we use Rosenbrock's function (a test-bench for global optima-seeking) as the cost function\;(See Figure \ref{fig:RF}). The Rosenbrock's function is a non-convex\footnote{The Hessian at all points of the function are not positive semi-definite \cite{rock}.} function which is known for its hard to find minima. This global minima exists at $(1,1)$ located within a large valley making the optimisation computationally hard. Rosenbrock's function is given as: 
\begin{equation}
    f(x,y)\, = \,(1-x)^2+100\,(y-x^2)^2.
    \label{6}
\end{equation}
\noindent
As mentioned earlier, the performance of an optimisation scheme is characterised by the spectral condition number of the Hessian of the cost function calculated at its minima. The spectral condition number of a matrix given as $\kappa = \frac{\lambda_{max}}{\lambda_{min}}$, where $\lambda$ signifies the eigenvalues \cite{strang2006linear}.\\[0.1cm] At its minima, the spectral condition number of the Hessian of Rosenbrock's function is calculated to be $\kappa = 2508$.\cite{beck} This indicates that the system is ill-conditioned and most first order gradient-based methods would require a large number of iterations to find its minima. In the continuous-time method, it has been found that it can only be solved by considering implicit methods which can handle such stiff systems. \cite{stiff}
\begin{figure}
\centering
\includegraphics[width=1\columnwidth]{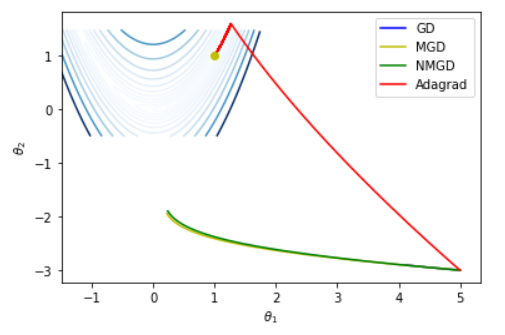}
 \caption{\small Some existing optimisation algorithms applied to Rosenbrock's function for the initial condition (5,-3).} 
 \label{fig:gg}
\end{figure}
\noindent
First order gradient algorithms when applied to Rosenbrock's function performs quite poorly in discrete-time; i.e. they are unable to find the minima within a practical run-time. While ADAM \cite{Kingma2015AdamAM} and other hybrid algorithms are able to solve this problem, they do so at a high computational cost; i.e. it takes them a large number of iterations to minimise Rosenbrock's function and often converge slowly due to the stiffness of the system as shown in Figure \ref{fig:gg} (note that ADAM takes at least $10^5$ iterations to converge to the minima for $\epsilon = 10^{-8}$ for Rosenbrock's function \cite{Study}).\\[0.1cm]
\noindent 
This paper is divided into two main sections to deal with the continuous and the discrete-time analyses of the proposed algorithm.
\footnote{Please note that in this paper, all simulations and results have been compiled in the MATLAB\texttrademark \,and SimuLink\texttrademark \,environments. All graphs where the X-axis has not been explicitly mentioned, denotes time in seconds.}
\noindent
\section{Continuous-Time Analysis} \label{sec:cont}
\subsection{Motivation}
\begin{figure}
\centering
    \includegraphics[width = 1\columnwidth]{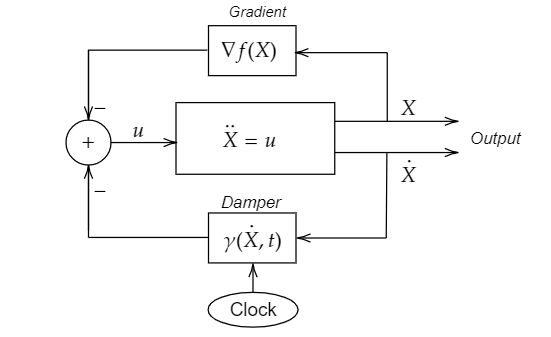}
\caption{Whiplash control system block diagram.} \label{fig:bkk}
\end{figure}

\noindent
The design of the proposed algorithm was inspired by a physical understanding of the dynamical system (\ref{7}) which involves control using its momentum \cite{10.5555/2051759}. This has also been motivated by the recent paper by Attouch et al. \cite{attouch2021fast}. In \cite{attouch2021fast}, closed-loop control is considered using multiple scenarios with a damping constant of the form $\gamma = r\,|\dot{x}|^{p-2}$, where $p$ and $r$ are positive constants (control parameters) and $\dot{x}$ is the velocity. However, the simulation results we obtained for this algorithm were not particularly satisfactory for the minimisation of Rosenbrock's function. We make two particular observations:
\begin{enumerate}
\item The damping function did not adequately stabilise the system over long intervals for various values of $p$ and $r$. This inspired us to use the control parameter $r$ as $t$ (time) and $p = 4$.\footnote{Note that this implies for the vector form $||\dot{X}||^2_2$ becomes the inner product of the velocity vector with itself.}
\item The lack of stability of the system over large intervals of time led us to look at (\ref{7}) as a linear time-variant system and to perform a corresponding pole placement.
\end{enumerate}
The linearised ODE for the method (\ref{6}) applied to Rosenbrock's function is equivalent to the form:
\begin{equation}\label{8}
\ddot{X}+\gamma\dot{X}+\nabla^2f\,X=0,
\end{equation}
where the Hessian $\nabla^2 f$ for Rosenbrock's function at the minima $X=X^*$, is:
\begin{equation}\label{9}
\nabla^2\,f =
\begin{bmatrix}
 802 & -400\\
 -400 & 200
\end{bmatrix}.
\end{equation}
Hence, we look at the eigenvalues for the linear time variant system matrix of (\ref{8}) to understand its convergence rate \cite{tomas2010linear}. For this purpose, we define a new system with state $y$ for the underlying linearised state space system $\dot{x} = A\,x$, which is given as:
\begin{equation}\label{10}
    y = x\,e^{-\eta t}, \;\; \eta>0\,.
\end{equation}
This leads to the redefined system:
\begin{equation}\label{11}
     \dot{y} + \eta\,y = A\,y\,,
\end{equation}
which in matrix form can be written as:
\begin{equation}\label{12}
  \dot{y} = (A\,-\eta I)\,y.
\end{equation}
The minimum real part of the eigenvalues of the system (11) can be written as:
\begin{equation}\label{13}
  \sigma = \frac{-\,\gamma}{2}-\eta
\end{equation}
which indicates an exponential convergence rate of at least $\eta$. Based on this, we choose the value of $\eta = 1$.
Thus we arrive at our hybrid gradient descent optimisation method which we shall refer to as the \textit{whiplash inertial gradient optimisation method} (\ref{14}):
\begin{equation}\label{14}
    \ddot{X}+(1+t\,||\dot{X}||^2\,)\dot{X}+\nabla f(X) = 0.
\end{equation}
A block diagram for this algorithm is shown in Figure \ref{fig:bkk}. The motivation for this nomenclature can be found in subsection \ref{sr}.

\subsection{Convergence Analysis using the Lyapunov Method}
\noindent
For an autonomous system, of the form:
\begin{equation}\label{15}
    \dot{x} = f(x) 
\end{equation}
we can guarantee global asymptotic stability if there exists a functional, $V\,(x)$, such that:
\begin{equation} \label{16}
\begin{split}
 &  V(x)>0\;\; \forall \,x\neq 0 \\
 &  V(0)=0 \\
 &  \dot{V}(x)\,<\,0\;\; \forall\,x \neq 0\,.
\end{split}
\end{equation}
Furthermore, we consider the La-Salle principle of invariance \cite{Khalil:1173048} and suppose there exists a continuously differentiable, positive definite, radially unbounded function  $V(z)\,: \mathrm{R}^n \rightarrow \mathrm{R}$ such that $\forall \;z \in \mathrm{R^n}$:
\begin{equation}\label{17}
  \frac{\partial{V}}{\partial{x}}(z-x_e)\,f(z)\leq W(z) \leq 0\,.
\end{equation}
Then, $x_e$ is a Lyapunov stable equilibrium point, and the solution always exists globally. \say{Moreover, $x(t)$ converges to the largest invariant set M contained in $E\,=\,\{ z \in \mathrm{R^n} \, :W(z)=0\}$. When $ W(z)=0$  only for $ z=x_{e}$ then $ E=\{x_e\}$. Since $ M \subset E$ therefore $x(t) \rightarrow x_e $ which implies asymptotic stability. Even when $E \neq \{x_e\}$, we often have the condition  $M=\left\{x_e\right\}$ from which we can conclude asymptotic stability} \cite{Kalman1960}.
This is used in our analysis for the general inertial gradient dynamical system \cite{ATTOUCH20175412} by defining a candidate Lyapunov function $W(t)$ for all damping functions $\gamma$ such that:
\begin{equation}\label{18}
W(t) = \frac{1}{2} ||\dot{x}||^2 + f(x) - f(x^*),
\end{equation}
where $f(x^*)$ denotes the minima of the cost function, which satisfies all the conditions in (\ref{11}). Upon replacing the time derivative in the equation (\ref{7}) for $\ddot{X}$, we obtain:
\begin{equation}\label{19}
   \begin{split}
    &  \dot{W}(t) = \langle \dot{X},\ddot{X} \rangle + \langle \nabla(f),\dot{X}
    \rangle \\
    & = - \; \gamma ||\dot{X}||^2\; \leq 0\ \; \forall \, \gamma \geq 0\,.
    \end{split} 
\end{equation}
This shows that the time derivative of our Lyapunov candidate is negative semi-definite. This shall suffice to show using La-Salle's principle of invariance that the set of accumulation points of any trajectory is contained in $\mathcal{I}$, where $\mathcal{I}$ is the union of complete trajectories contained entirely in the set $\{\mathbf{x} :\dot {W}(\mathbf{x})=0 \} $\cite{Khalil:1173048}. Thus by Lyapunov's Second Theorem, we have the functional $W$ is positive definite; i.e. \say{$\mathcal{I}$ contains no trajectory of the system except the trivial trajectory $\mathbf{x}(t) \equiv \mathbf{0}$ and as $W$ is radially unbounded; i.e. $W(\mathbf{x} )\rightarrow \infty$ as $\mathbf{||x||} \rightarrow \infty$}, we conclude that the origin is \textit{globally asymptotically stable} \cite{Kalman1960}.
\subsection{Simulation Results}\label{sr}
\begin{figure}\centering
\includegraphics[width=1\columnwidth]{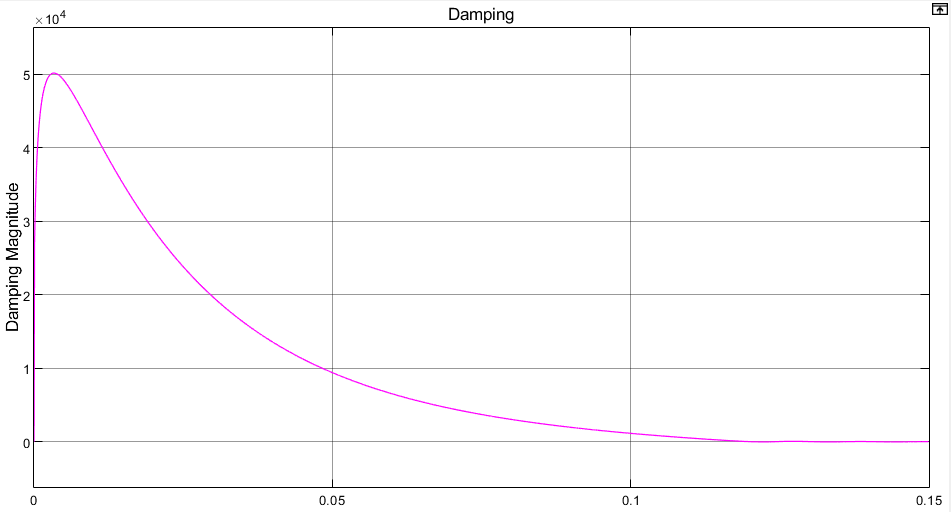}
 \caption{\small The Closed-Loop Damping coefficient $\gamma$.}
 \label{fig:5}
\end{figure}

\begin{figure}
    \centering
    \subfigure[]{\includegraphics[width=0.35\textwidth, height=0.25\textwidth]{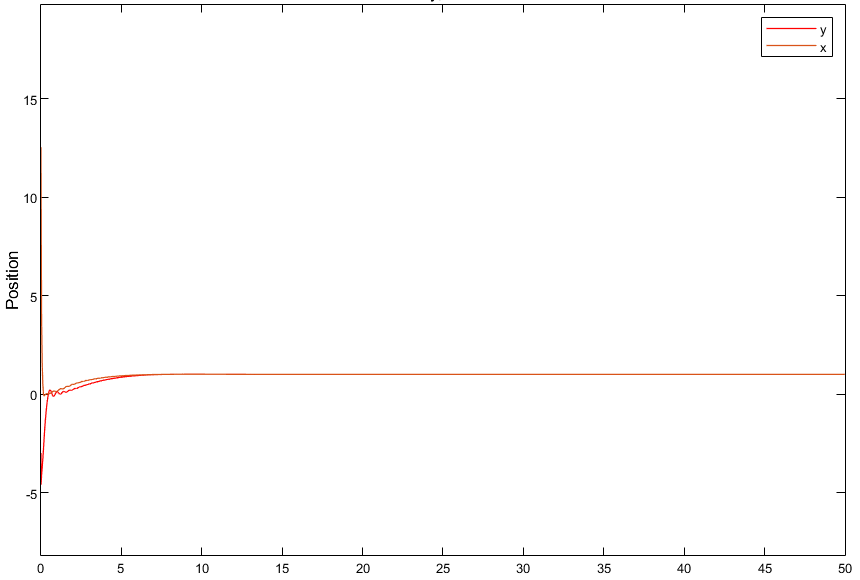}} \\[0.1cm]
    \subfigure[]{\includegraphics[width=0.35\textwidth, height=0.25\textwidth]{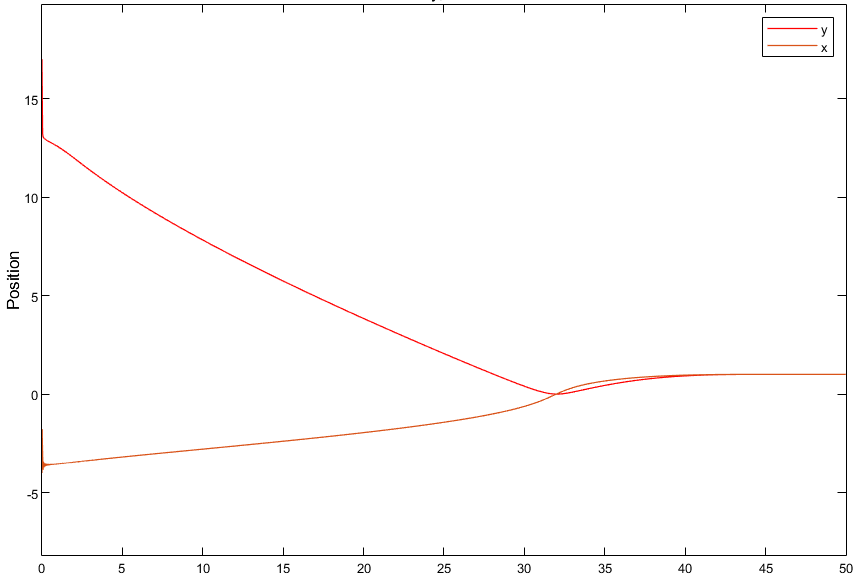}} \\[0.1cm]
    \caption{\small  Simulations for the continuous-time whiplash gradient descent: (a) Initial condition (12,-3) (b) Initial condition (-4,17)}
    \label{fig:6}
\end{figure}

\noindent
We modelled the whiplash inertial gradient dynamic system (14) using the Euler fixed step integrator (\texttt{ode1}) on Rosenbrock's function, using a step-size of $0.001$. A few examples of the state trajectories have been shown in Figure \ref{fig:6} starting from different initial conditions. While all starting speeds achieve convergence (unlike the Nesterov scheme which must be started at $\dot{X}(0) = 0$ \cite{ATTOUCH20175412}\cite{su2015differential} to ensure convergence) it should be noted that we start the system at an arbitrary fixed velocity of $\dot{X}(0) = -1000$ which shows rapid convergence for all initial conditions\footnote{Further research will be required to study the effect of various starting speeds for this system.}.
\\[0.1cm]
If one looks at Figure \ref{fig:5}, where we have analysed the system's damping coefficient $\gamma = 1+t\,\langle \dot{X},\dot{X}\rangle$, it shows a sharp rise followed by an abrupt fall for high starting speeds (which ensures faster convergence). This phenomenon replicates the physical process of the \textit{whiplash} effect, and hence motivated its nomenclature. We hypothesise that (as explained by Kovachki and Stuart in their paper \cite{kovachki2021continuous}) this effect might be responsible for the rapid stabilisation of the system in the transient phase.
\section{Discrete-Time Analysis} \label{sec:disc}
\subsection{Discretisation}
\noindent
The discretization that is used is the semi-implicit or symplectic Euler method \cite{bp}. Using a discrete-time step $s$ and sampling of $t \approx \, k\sqrt s$, we obtain a two-state estimate of the acceleration and velocity as shown below:
\begin{equation}\label{20}
\begin{split}
& \dot{X} \approx v_k = \frac{x_{k+1}-x_{k}}{\sqrt{s}},\\[0.1cm]
& \ddot{X} \approx \frac{v_{k+1}-v_{k}}{\sqrt{s}},\\
& 1+t\,||\dot{X}||^2\, \approx 1+k\sqrt{s}||v_k||^2.
\end{split}
\end{equation}
Now, we modify (\ref{14}) to add a fixed mass to the system, which up until this point has been considered to be of unit magnitude. The choice of mass that we shall make is $m = \frac{1}{\sqrt{s}}$. This idea of introducing this mass in inertial gradient flow methods while discretising them has been inspired from the idea of selective mass scaling in finite element methods \cite{fem}, where the iterative process can be scaled by choosing an effective mass. The rationale behind this is since the discrete-time method depend heavily on the step-size, they take much longer to attenuate for smaller step-sizes. Hence, to counter this effect, we may introduce such a fixed mass, which scales the dynamics, depending on the step-size. Upon making these substitutions and modification to (\ref{14}), we obtain
\begin{equation}\label{21}
    (\frac{1}{\sqrt{s}})\frac{v_{k+1}-v_{k}}{\sqrt{s}} + (1 + k\sqrt{s}||v_k||^2)\frac{v_k}{\sqrt{s}} +\nabla f(x_k) = 0.
\end{equation}
We can re-write (\ref{21}) using (\ref{20}) as:
\begin{equation}
\label{22}
v_{k+1} = (1-\sqrt s - ks\;v_k^T \cdot v_k)v_k - s\nabla{f(x_{k})}\,.
\end{equation}
Now, we consider the symplectic approximation for the Lyapunov stable system \cite{hair} such that $\|\dot{x}(t)\|\rightarrow0$ as $t\rightarrow \infty$. This implies that $\|v_k\|\rightarrow0$ as $k\rightarrow\infty$. Therefore, we introduce $z_k = x_{k}-x_{k-1} = \sqrt{s}v_{k-1}$. For all asymptotic analyses, there is no practical difference between the sequences $z_k$ and $ v_k$ as
\begin{equation}
    \lim_{k\rightarrow \infty} \|z_{k+1}-v_k\| = |\sqrt{s}-1|\|v_k\| = 0.
\end{equation}
We may consider this as two transforms. First as a scaling of the system, followed by a backward recursion:
\begin{equation}
    \lim_{t\rightarrow\infty} |\sqrt{s}-1| \dot{x}(t)
\approx \max_{0<k\leq \frac{T}{\sqrt{s}}} (x(k\sqrt{s})-x((k-1)\sqrt{s})).
\end{equation}
This trick simplifies our system's updates while keeping intact the geometry of the dynamical system and does not change the global nature of the system's convergence\footnote{We have verified this claim using numerical results.}. This particular choice of design for the algorithm simplifies the computation and makes discrete-time analyses of convergence considerably easier. We finally have the consolidated scheme as:
\begin{equation}
\label{WGDA}
    \begin{cases}
    & x_1 = x_0 - s \nabla{f(x_0)},\\
    & z_k = x_{k}-x_{k-1},\\
    & \alpha_k = 1-\sqrt s - ks||z_k||^2, \\
    & z_{k+1} = \alpha_k\,z_k - s\nabla f(x_{k}).
    \end{cases}
\end{equation}
\subsection{Algorithm}
\noindent
This discrete-time scheme (\ref{WGDA}) can be translated to the following algorithm\footnote{The code is available on the licensed repository: \url{https://github.com/SubhransuSekharBhattacharjee-01/Whiplash.git}} using a step-size $s$ and $n$ iterations, with initial starting point $x_0$ and final point $x_n$ ($ \ref{alg:1} $). Unlike prior gradient descent algorithms, which are capable of minimising Rosenbrock's function, this algorithm does not use any hyper-parameters. Instead, it uses a simple two-step assignment to update the discrete-time damping in every iteration. The \textit{zeroth} step (\ref{A1}) of the iteration $x_1$ is a gradient descent step \cite{beck} which assigns the initial momentum for the first iteration as $z_1 = x_1 -x_0$.
 \begin{algorithm}
 \caption{The whiplash gradient descent algorithm}
\textbf{Input: $\nabla{f(x)}, n, s, x_0$}
 \label{alg:1}
 \begin{algorithmic}[1]
    \STATE \label{A1} \textbf{Initialise: $x_1 \leftarrow x_0 - s \nabla{f(x_0)}$}
    \STATE $k = 1$
    \WHILE{$k \leq n$}
    \STATE $z_k \leftarrow x_{k}-x_{k-1}$
    \STATE $\alpha_k \leftarrow 1- \sqrt s - ks(z_k^T \cdot z_k)$
    \STATE $x_{k+1} \leftarrow x_{k} +\alpha_k z_k - s\nabla{f(x_{k}) }$
    \STATE $(x_{k-1},x_{k}) \leftarrow (x_{k},x_{k+1})$
    \STATE $k \leftarrow k+1$
    \ENDWHILE\\
     \end{algorithmic}
\textbf{Output: $x_n$} 
 \end{algorithm}
\subsection{Numerical Results}
\begin{figure}
    \centering
    \includegraphics[width = 1\columnwidth]{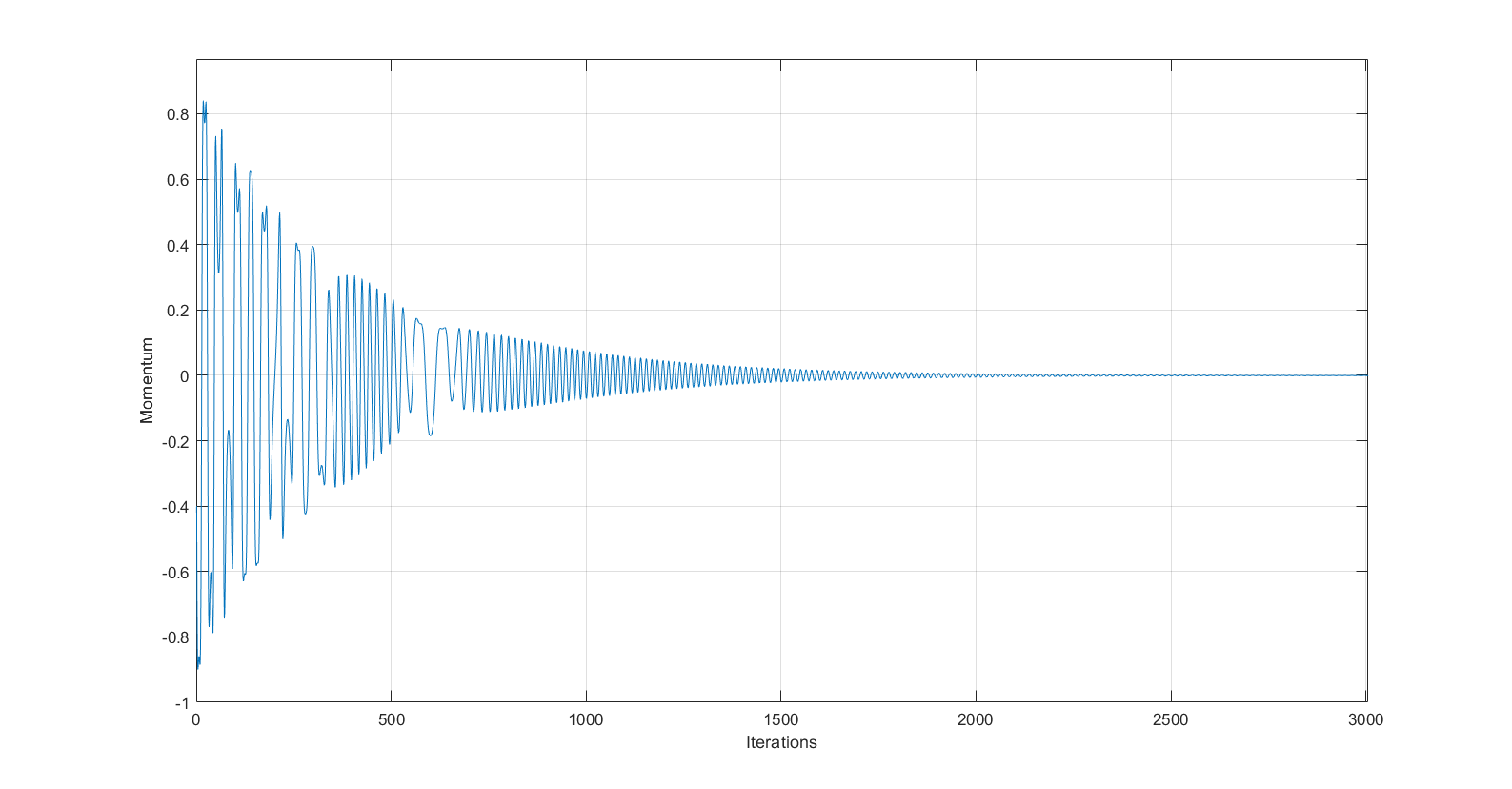}
    \caption[]
    {\small Momentum plot: Convergence of the momentum with iterations.}
    \label{fig:8}
\end{figure}

\begin{figure*}
        \centering
            \includegraphics[width = 2\columnwidth]{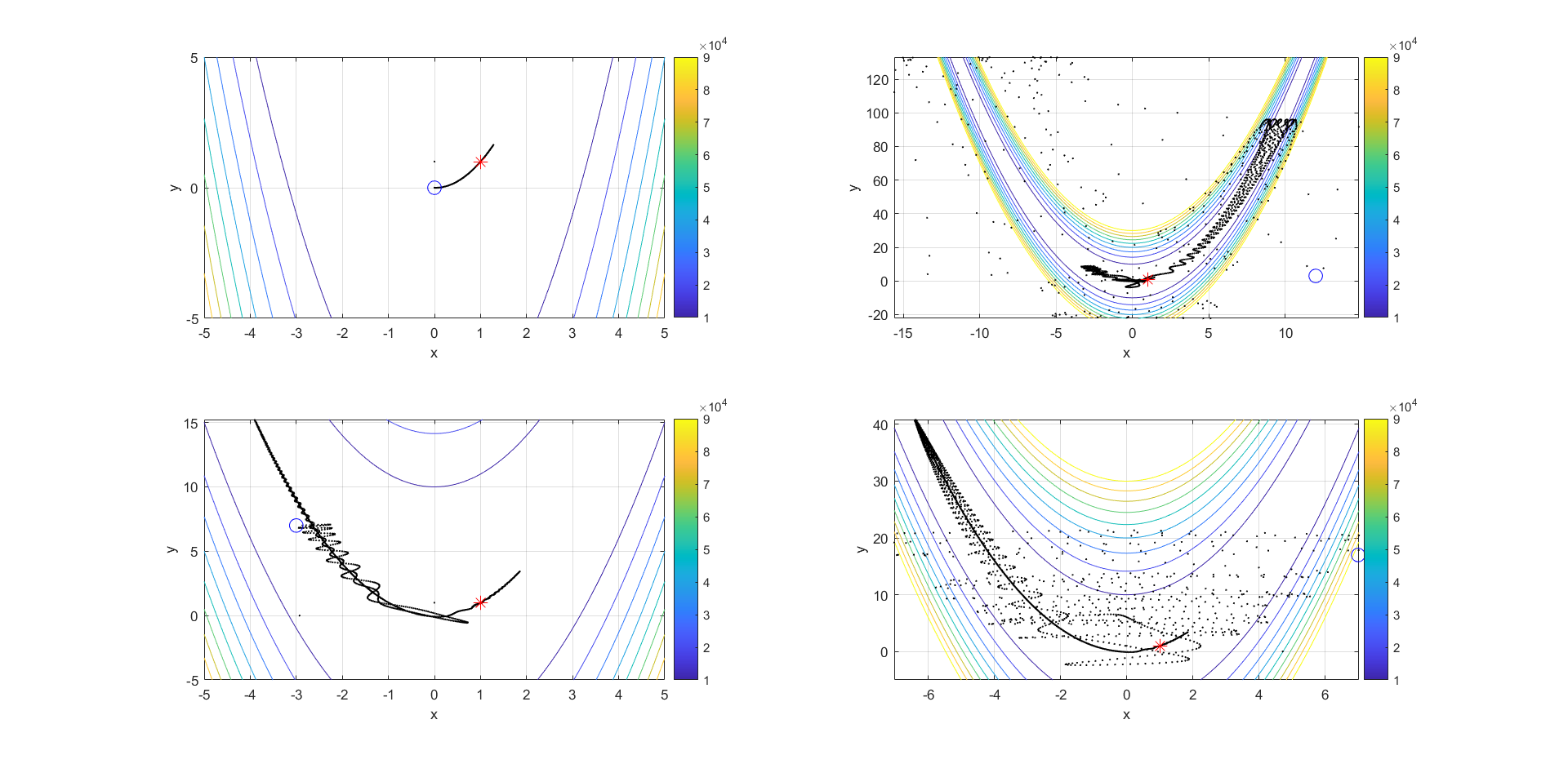}
        \caption[]
        {\small  Trajectories dependent on initial condition: Note the blue circle denotes the starting point and the red star denotes the ending point.
        From upper left corner the initial conditions are (0,0), upper right(12,3), lower left (-3,7) and lower right (-7,17).}
        \label{fig:7}
    \end{figure*}

\begin{figure*}
    \centering
    \includegraphics[width = 1.35\columnwidth]{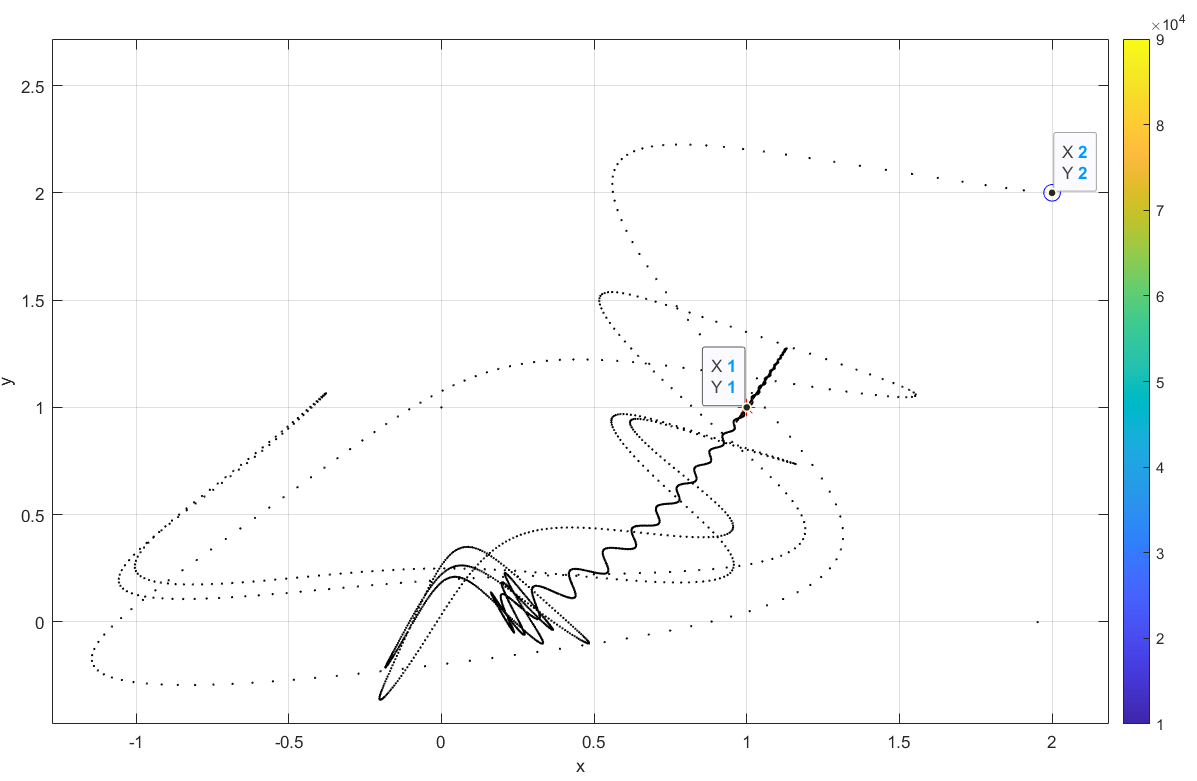}
    \caption[]
    {\small Zoomed Trajectory Plot with starting point (2,2).}
    \label{fig:9}
\end{figure*}

\noindent
Numerical results for the whiplash gradient descent algorithm applied to Rosenbrock's function are promising. As explained previously, for the stiffness of the system to be taken into account, we need a step-size of no more than $10^{-5}$. This is because the algorithm is unable to learn the gradient of the cost function and picks up momentum without correcting the damping. For a sufficiently small step-size, the whiplash gradient descent algorithm successfully found the minima of Rosenbrock's function for all initial conditions.\\[0.1cm]
We have shown a few examples in Figure \ref{fig:7}. A plot of the momentum growth Figure \ref{fig:8} shows a saturation effect. This indicates that Rosenbrock's function optimisation has been achieved over the given time interval. Figure \ref{fig:9} shows the nature of the trajectory, as it approaches the minima.
\section{Conclusion and Future Work} \label{sec:Conclusion}
\noindent
From the above results, we can see that the proposed whiplash gradient descent algorithm is capable of fast optimisation of Rosenbrock's function. We constructed this algorithm using a non-linear controller motivated by the nature of the momentum-control structure. However, it must be understood that this controller might not be optimal even for Rosenbrock's function. This is because unlike linear systems, where we could predict the results from clear theoretical motivations, we do not have any such tools for analysis for the non-linear case.\\[0.1cm] Thus, as a direction for future research, we will need to reconsider the classical Lur'e problem, for the absolute stability of the entire class of the inertial gradient systems, involving a feedback path that contains a sector-bound non-linearity \cite{haddad2011nonlinear}. We will also need to research further to understand the theoretical and practical limitations of closed loop control for the generalised inertial gradient system. Furthermore, we need to study the effect of variation in starting speed and the hypothesis regarding the stabilisation effect. Deriving upper bounds for the rates of convergence, using Lyapunov arguments, will be another direction for research. For that we may consider a Lyapunov argument of the form:
\begin{equation}\label{23}
    \mathscr{E}_t = \mathscr{P}\,(f-f^*) + \frac{1}{2}||x-x^*+\mathscr{L}||^2,
\end{equation}
where $f^*$ denotes the optimal value of the cost function where $\{x^*:f^*=f(x^*)\}$, $\mathscr{P}\,(t,\alpha)$ is the \textit{converger to the scheme} \cite{Baes}, with the hyper-parameter $\alpha$ and $\mathscr{L}(x,\dot{x},t)$ is a differentiable non-linear function. Such a Lyapunov candidate might be used to estimate the convergence rate whose upper bound is $\mathcal{O}(\frac{1}{\mathscr{P}})$, for various geometrical considerations like the nature of convexity of the cost function, only if $\dot{\mathscr{E}_t}<0$ \cite{awil}.
\section{Acknowledgement}
\noindent
We thank the College of Engineering and Computer Science, Australian National University for supporting our work.
\bibliographystyle{ieeetran}
\bibliography{biblio.bib}
\end{document}